\documentclass[12pt, draft, leqno]{amsart}
\usepackage[mathscr]{eucal}
\usepackage{amssymb}

\renewcommand{\baselinestretch}{1.5}

\theoremstyle{plain}
\newtheorem{theorem}{Theorem}
\newtheorem{corollary}[theorem]{Corollary}

\newtheorem{proposition}[theorem]{Proposition}

\theoremstyle{remark}

\newcommand\rng{\operatorname{rng}}
\newcommand\tr{\operatorname{tr}}
\newcommand\R{\mathbb R}

\begin{document}
\vglue -20pt

\centerline{\large{\textbf{A CLASS OF PRESERVERS ON HILBERT SPACE
EFFECTS}}}

\centerline{\large{\textbf{INCLUDING ORTHO-ORDER AUTOMORPHISMS}}}

\centerline{\large{\textbf{AND SEQUENTIAL AUTOMORPHISMS}}}

\vskip 10pt
\centerline{LAJOS MOLN\'AR}

\vskip -5pt
\centerline{Institute of Mathematics and Informatics}

\vskip -5pt
\centerline{University of Debrecen}

\vskip -5pt
\centerline{4010 Debrecen, P.O.Box 12}

\vskip -5pt
\centerline{Hungary}

\vskip 5pt
\centerline{e-mail: \texttt{molnarl@math.klte.hu}}

\vskip 10pt
\centerline{Running title:}

\centerline{\scshape{PRESERVERS ON HILBERT SPACE EFFECTS}}

\pagestyle{myheadings}
\markboth{\textsc{\SMALL PRESERVERS ON HILBERT SPACE EFFECTS}}
{\textsc{\SMALL PRESERVERS ON HILBERT SPACE EFFECTS}}

\normalsize
\vskip 10pt
\centerline{\textsc{Abstract}}
In this paper we study a new class of transformations on the set of all
Hilbert space effects. This consists of the bijective maps which
preserve the order and zero product in both directions. The main
result of the paper gives a complete description of the structure
of those transformations. As applications we obtain additional new
results and some former ones as easy corollaries. In particular,
we obtain the form of the
ortho-order automorphisms as well as that of the sequential
automorphisms. In the last
paragraph of the paper we show that these two kinds of
automorphisms belong to our class of transformations even when
their domain is the set of all effects in a general von
Neumann algebra.
\noindent


\newpage
\renewcommand{\baselinestretch}{1.5}
\normalsize

\section{Introduction and Statements of the Results}

Effects play very important role in the quantum theory of
measurement (see, for example, \cite{BusLahMit91}). In the Hilbert
space formalism of the theory, the so-called Hilbert space effects
are the positive bounded linear operators on a Hilbert space $H$
which are bounded by the identity $I$. The set of all effects on
$H$ is denoted by $E(H)$.

There are several operations and relations defined on $E(H)$ which
are important from different aspects of the theory. What concerns
the present paper, here we are interested in the ortho-order
structure and in the sequential structure on $E(H)$. The first one
is obtained as follows. The usual ordering $\leq$ among
self-adjoint bounded linear operators gives rise to a partial
order on $E(H)$ and the operation $':A \mapsto I-A$ defines a kind
of orthocomplementation on $E(H)$. This relation and operation
together give the ortho-order structure on $E(H)$ \cite{Lud83}. As
for the second structure, it comes from the sequential product
which is defined as follows. If $A,B\in E(H)$, then their
sequential product is $A\circ B=\sqrt{A}B\sqrt{A}$
\cite{GudNag01}.

Supposing $\dim H\geq 3$, the automorphisms of $E(H)$ with respect
to the ortho-order structure (called ortho-order automorphisms) as
well as the ones with respect to the sequential product (called
sequential automorphisms) are known to be implemented by unitary
or antiunitary operators. This means that all those automorphisms
are of the form
\begin{equation}\label{E:hat10}
\phi(A)=UAU^* \qquad (A\in E(H))
\end{equation}
where $U$ is a unitary or antiunitary operator. The result concerning
ortho-order automorphisms was obtained by Ludwig in
\cite[Section V.5]{Lud83} (the proof was later clarified in
\cite{CasVitLahLev00}) while the corresponding result on
sequential automorphisms was given by Gudder and Greechie in
\cite{GudGre02a}. In the paper \cite{ML01c} we showed that
Ludwig's theorem holds also in the 2-dimensional case.

In our recent paper \cite{ML01d} we presented some characterizations of
the ortho-order automorphisms of $E(H)$ by means of their
preserver properties. This investigation was motivated by the
extensive study of preserver problems in matrix theory and in
operator theory. Roughly speaking, preserver problems are
concerned with the description of all
transformations (called preservers) on a
given set (preferably equipped with an algebraic structure) which
preserve a certain relation between the elements, or a given
subset of elements, or a quantity attached to the elements of the
underlying set.
In the present paper we continue the investigation started in \cite{ML01d}
by studying the bijective maps $\phi :E(H) \to E(H)$ which preserve the order and
zero product in both directions, i.e., which satisfy
\begin{equation}\label{E:hat20}
A\leq B \Longleftrightarrow \phi(A)\leq \phi(B) \qquad
(A,B\in E(H))  \tag{a}
\end{equation}
and
\begin{equation}\label{E:hat21}
AB=0 \Longleftrightarrow \phi(A)\phi(B)=0 \qquad
(A,B\in E(H)). \tag{b}
\end{equation}
In what follows we present a complete description of the
structure of those maps and give several
applications of the corresponding result. Among others, we obtain
the form of all bijective transformations of $E(H)$ which
\begin{itemize}
\item[(i)] preserve the order in both directions and
\item[(ii)] map one single nontrivial scalar operator
$\lambda I$ $(\lambda \neq 0,1)$ to an operator of the same type.
\end{itemize}
At the first glance, it might be rather surprising that such maps are of a nice form
but it turns out that they in fact satisfy
\eqref{E:hat20} and \eqref{E:hat21}. Next, we easily recover
one of the main results in \cite{ML01d} on maps preserving
the order and coexistency in both directions (the definition of coexistency
is given below). What is probably
more important, we also show that the main result of the present paper readily implies the
former results on the structure of ortho-order automorphisms
and sequential automorphisms of $E(H)$ given in
\cite{Lud83} and in \cite{GudGre02a}, respectively. In fact, we prove
that those automorphisms also satisfy \eqref{E:hat20} and \eqref{E:hat21}.
It should be emphasized that, due to the proof of our main result, all
statements presented in the paper are valid in
2-dimensional case as well. In particular, this holds for the result on the
form of sequential automorphisms of $E(H)$ which is a
new result.

Finally, in the last paragraph of the paper we demonstrate that even on
the set $E(\mathcal A)$ of all effects belonging to a general von
Neumann algebra, the ortho-order automorphisms and the
sequential automorphisms (the definitions should be
self-explanatory) belong to our new class of preservers, i.e.,
they preserve the order and zero product in both directions. This
observation is worth mentioning since in that generality it is
quite hard to see any connection between those two kinds of
automorphisms. Therefore, we believe that our preservers deserve
attention and it has sense to study them in other contexts as
well.

Now we turn to the precise formulations of the results of the paper.
\emph{In what follows we
assume that $H$ is a (complex) Hilbert space with $\dim H\geq 2$.}
Our main result reads as follows.

\begin{theorem}\label{T:hat2}
Let $\phi:E(H)\to E(H)$ be a bijective map which preserves the
order and zero product in both directions. Then there is an either
unitary or antiunitary operator $U$ on $H$ and a real number $p <1$
such that with the function $f_p (x)=\frac{x}{xp +(1-p)}$ $(x\in
[0,1])$ we have
\[
\phi(A)=Uf_p (A) U^* \qquad (A\in E(H)).
\]
Here, $f_p (A)$ denotes the image of the function $f_p$ under the
continuous function calculus belonging to the operator $A$.
\end{theorem}

To get the form of all bijective maps on $E(H)$ with the
properties (i), (ii), we need the following proposition which
might be interesting on its own right.

\begin{proposition}\label{P:hat1}
Let $\phi:E(H)\to E(H)$ be a bijective map which preserves the
order in both directions and suppose that there is a single pair
of scalars $\lambda, \mu \in ]0,1[$ such that $\phi(\lambda I)=\mu
I$. Then $\phi$ preserves zero product in both directions.
\end{proposition}

The above results have the following immediate consequences. To
the second statement in the corollary below observe that scalar
operators in $E(H)$ can be characterized as effects which commute
with every other effect.

\begin{corollary}\label{C:hat3}
Let $\phi:E(H)\to E(H)$ be a bijective map which preserves the
order in both directions and suppose that there is a single pair
of scalars $\lambda, \mu \in ]0,1[$ such that $\phi(\lambda I)=\mu
I$. Then there exist an either unitary or antiunitary operator $U$ on
$H$ and a real number $p <1$ such that with the function $f_p
(x)=\frac{x}{xp +(1-p)}$ $(x\in [0,1])$ we have
\[
\phi(A)=Uf_p (A) U^* \qquad (A\in E(H)).
\]
In particular, we obtain the same form for any bijection of $E(H)$
which preserves the order and commutativity in both directions.
\end{corollary}

It is easy to see that if a function $f_p$ appearing above satisfies
$f_p(\lambda)=\lambda $ for some $\lambda \in ]0,1[$, then $p=0$
and hence $f_p$ is the identity. This gives us the following
corollary stating that if an order preserving bijection $\phi$ of
$E(H)$ fixes one single nontrivial scalar operator, then $\phi$ is
implemented by an either unitary or antiunitary operator.

\begin{corollary}\label{C:hat4}
Let $\phi:E(H)\to E(H)$ be a bijective map which preserves the
order in both directions and suppose that there is a scalar
$\lambda \in ]0,1[$ such that $\phi(\lambda I)=\lambda I$. Then
there exists an either unitary or antiunitary operator $U$ on $H$ such
that
\[
\phi(A)=U A U^* \qquad (A\in E(H)).
\]
\end{corollary}

The following corollary in the case when $\dim H\geq 3$ appears in
\cite{ML01d} as Theorem 1. Here we shall present a short proof of it
based on our main result which works also in the 2-dimensional
case. Recall that two effects $A,B\in E(H)$ are called coexistent
if they are in the range of a positive operator valued measure or,
equivalently,
if there are effects $E,F,G\in E(H)$ such that
\[
A=E+G, \quad B=F+G, \quad E+F+G \in E(H).
\]

\begin{corollary}\label{C:hat5}
Let $\phi:E(H)\to E(H)$ be a bijective map which preserves the
order and coexistency in both directions. Then there is an either
unitary or antiunitary operator $U$ on $H$ such that
\[
\phi(A)=U A U^* \qquad (A\in E(H)).
\]
\end{corollary}

As for the last two corollaries below, we shall see that the structure
of the ortho-order automorphisms and that of the sequential
automorphisms of $E(H)$ can be deduced from our main result even in the
2-dimensional case. The proofs are based on the observation that
both kinds of automorphisms preserve the order and zero product in
both directions. (For a similar observation concerning effects in
general von Neumann algebras, see the last paragraph of the
paper.) Recall that a bijective map $\phi: E(H) \to E(H)$ is
called an ortho-order automorphism if for any $A,B\in E(H)$ we
have
\[
A\leq B \Longleftrightarrow \phi(A)\leq \phi(B),
\]
(i.e., $\phi$ preserves the order in both directions) and
\[
\phi(A')=\phi(A)'
\]
(i.e., $\phi$ preserves the orthocomplements). Moreover, a bijective
map $\phi: E(H) \to E(H)$ is called a sequential automorphism if
\[
\phi(A\circ B)=\phi(A)\circ \phi(B)
\]
holds for every $A,B\in E(H)$.

\begin{corollary}\label{C:hat6}
Let $\phi:E(H)\to E(H)$ be an ortho-order automorphism. Then
$\phi$ preserves the order and zero product in both directions. In
fact, there is an either unitary or antiunitary operator $U$ on
$H$ such that
\[
\phi(A)=U A U^* \qquad (A\in E(H)).
\]
\end{corollary}

\begin{corollary}\label{C:hat7}
Let $\phi:E(H)\to E(H)$ be a sequential automorphism. Then $\phi$
preserves the order and zero product in both directions. In fact,
there is an either unitary or antiunitary operator $U$ on $H$ such
that
\[
\phi(A)=U A U^* \qquad (A\in E(H)).
\]
\end{corollary}

Finally, we remark that the converse statements in all of the
above results excluding Proposition~\ref{P:hat1} are also true. We
mean that if a transformation is of the form which appears in the
formulation of the corresponding result, then it has the
properties which were required there. Since the verification of
this observation needs only elementary computations, hence we omit
them.

\section{Proofs}

First we emphasize that the proof of our main
result has many common points with the proof of
\cite[Theorem]{ML01c}. In fact, the argument to be presented below
can be considered as an
adaptation of the proof given in \cite{ML01c} for another situation. So,
it would have been a possibility just to point out where and how the
proof in \cite{ML01c} should be modified to get the statement of
the main result of this paper but we think that such a proof is quite
hard to follow and hence it does not meet the most elementary requirements.
Hence, we decided to present a complete, self-contained proof.

We begin with some notation and useful facts that we shall apply
in our arguments.

First, we note that the concept of the strength of an effect along
a ray plays very important role in what follows. This notion was
introduced by Busch and Gudder in \cite{BusGud99}.  If $A$ is an
effect on $H$, $\varphi$ is a unit vector in $H$ and $P_\varphi$
is the rank-one projection onto the subspace generated by
$\varphi$, then the quantity
\[
\lambda(A, P_\varphi)= \sup \{ \lambda \in [0,1] \, : \, \lambda
P_\varphi \leq A\}
\]
is called the strength of $A$ along the ray represented by
$\varphi$. Due to \cite[Theorem 4]{BusGud99} there is a very
useful formula to compute the strength. In fact, we have
\begin{equation}\label{E:hat3}
\lambda(A, P_\varphi)=
\left\{%
\begin{array}{ll}
    \|A^{-1/2}\varphi\|^{-2}, & \hbox{{\text{if }} $\varphi \in \rng (A^{1/2})$;} \\
    0, & \hbox{{\text{else}}.} \\
\end{array}%
\right.
\end{equation}
(Here, $\rng$ denotes the range of an operator and $A^{-1/2}$
denotes the inverse of $A^{1/2}$ on its range.)

Let $\phi:E(H)\to E(H)$ be a bijective map which preserves the
order in both directions. It was proved by Ludwig in \cite[Theorem
5.8., p. 219]{Lud83} that $\phi$ necessarily preserves the
projections in both directions. It is then trivial to see
that $\phi$ also preserves the rank of the projections (cf. the
proof of \cite[Theorem 1]{ML01d}).

An easy fact follows what we shall use several times. Namely,
if $A,B$ are effects, $B$ is of rank one and $A\leq B$, then $A$
is a scalar multiple of $B$. This observation and
the previous one have, among others, the following corollary.
Let $\phi$ be as
above, i.e., suppose that it is a bijection of $E(H)$ which
preserves the order in both directions. Then for every rank-one
projection $P$, there is a function $f_P :[0,1] \to [0,1]$ such
that
\[
\phi(tP)=f_P(t) \phi(P) \qquad (t\in [0,1]).
\]
By the order preserving property of $\phi$ and $\phi^{-1}$ we see
that $f_P$ is strictly increasing and bijective. In fact, we
have
\begin{equation}\label{E:hat5}
\phi^{-1}(t \phi(P))=f_P^{-1}(t)P \qquad (t\in [0,1]).
\end{equation}

Now, we turn to the proofs. In the proof of our main result
Theorem~\ref{T:hat2} we need the following proposition which
presents the solution of a functional equation.

\begin{proposition}\label{P:hat8}
Let $f, g:]0,1[ \to ]0,1[$ be functions and suppose that $f$ is a
strictly monotone increasing bijection. Let
\begin{equation}\label{E:FE}
f\biggl( \frac{x}{x+(1-x)y}\biggr)= \frac{f(x)}{f(x)+(1-f(x))g(y)}
\qquad (x,y \in ]0,1[).
\end{equation}
Then there are positive real numbers $a,b,c$ such that
\[
  f(x)=\frac{x^c}{x^c+a(1-x)^c} \qquad(x\in]0,1[)
\]
and
\[
  g(y)=b y^c\qquad(y\in]0,1[).
\]
\end{proposition}

\begin{proof} First note that since the function $f$ is continuous, equation \eqref{E:FE}
implies the continuity of $g$.

Next observe that with the notation
\[
  \begin{array}{rcll}
    \alpha(t)&=&\displaystyle{\frac{1}{1+e^t}} \qquad &(t\in \R), \\[3mm]
    \beta(x)&=&\ln\displaystyle{\frac{1-x}{x}} &(x\in]0,1[), \\[3mm]
    \gamma(y)&=&\ln y &(y\in]0,1[),
  \end{array}
\]

\noindent we have the identity
\[
  \frac{x}{x+(1-x)y}
   =\frac{1}{1+\exp\bigl(\ln\frac{1-x}{x}+\ln y\bigr)}
   =\alpha\bigl(\beta(x)+\gamma(y)\bigr)
\]
for all $x,y\in]0,1[$. Therefore, equation \eqref{E:FE} can be
rewritten as
\begin{equation}\label{E:FE2}
  f\circ\alpha\bigl(\beta(x)+\gamma(y)\bigr)
      =\alpha\bigl(\beta\circ f(x)+\gamma\circ g(y)\bigr)
      \qquad (x,y\in]0,1[).
\end{equation}
Substituting $x=\beta^{-1}(u)$ and $y=\gamma^{-1}(v)$ into
\eqref{E:FE2} and applying the inverse function of $\alpha$ to
both sides of \eqref{E:FE2}, we get
\begin{equation}\label{E:FE3}
  \alpha^{-1}\circ f\circ\alpha(u+v)
      =\beta\circ f\circ \beta^{-1}(u)
             +\gamma\circ g\circ\gamma^{-1}(v)
\end{equation}
for all $u\in\R$ and $v\in]-\infty,0[$. This means that the functions
\[
  F=\alpha^{-1}\circ f\circ\alpha,\quad
  G=\beta\circ f\circ \beta^{-1},\quad\mbox{and}\quad
  H=\gamma\circ g\circ\gamma^{-1}
\]
satisfy the following so-called Pexider equation
\[
  F(u+v)=G(u)+H(v) \qquad(u\in\R,\, v\in]-\infty,0[).
\]
Then, by known results of the theory of functional equations (cf.
\cite{Acz66}, or \cite{Kuc85}) and by the continuity of $F$, $G$,
$H$, it follows that there exist constants $a,b,c\in\R$ such that
\begin{eqnarray}
  F(w)&=&cw+a+b \qquad(w\in\R),\nonumber\\
  G(u)&=&cu+a\qquad\qquad(u\in\R),\label{E:FE5}\\
  H(v)&=&cv+b\qquad\qquad(v\in]-\infty,0[).\label{E:FE6}
\end{eqnarray}
Using \eqref{E:FE5} and the definition of $G$, we get that
$\beta\circ f(x)=c\beta(x)+a$. Easy computation yields that
\[
  f(x)=\frac{x^c}{x^c+e^{a}(1-x)^c} \qquad(x\in]0,1[).
\]
Similarly, the definition of $H$ and $\gamma$, and equation
\eqref{E:FE6} give
\[
  g(y)=e^b y^c\qquad(y\in]0,1[).
\]
The function $f$ being strictly increasing, $G$ is also increasing
and hence we get $c>0$.
\end{proof}

Now, we are in a position to prove our main result.

\begin{proof}[Proof of Theorem~\ref{T:hat2}]
The clue of the proof is to show that the functions $f_P$ (see the
first part of this section) do not depend on the rank-one
projections $P$. This will be done in what follows.

Let $P, Q$ be arbitrary mutually orthogonal rank-one projections.
By the order, rank and orthogonality preserving properties of
$\phi$ on set of all projections we clearly have
\[
\phi(P+Q)=\phi(P)+\phi(Q).
\]
Let $\lambda \in [0,1]$. From the inequality
\[
\phi(Q)\leq \phi(\lambda P+Q)\leq \phi(P+Q)=\phi(P)+\phi(Q)
\]
we obtain
\[
0\leq \phi(\lambda P+Q)-\phi(Q)\leq \phi(P).
\]
As $\phi(P)$ is of rank-one, according to the introduction of the
present section, this implies that there is a scalar $h_P(\lambda)
\in [0,1]$ such that
\[
\phi(\lambda P+Q)-\phi(Q)=h_P(\lambda) \phi(P)
\]
or, equivalently, that
\[
\phi(\lambda P+Q)= h_P(\lambda) \phi(P) +\phi(Q).
\]
Since $\phi^{-1}$ has the same properties as $\phi$, it can be
seen that the function $h_P:[0,1] \to [0,1]$ is a strictly
monotone increasing bijection. In fact, we have
\begin{equation}\label{E:hat6}
\phi^{-1}(\lambda\phi(P)+\phi(Q))=h_P^{-1}(\lambda)P+Q.
\end{equation}
We assert that $h_P=f_P$. Indeed, since
\[
f_P(\lambda)\phi(P)=\phi(\lambda P)\leq \phi(\lambda P+Q)=h_P(
\lambda) \phi(P)+\phi(Q),
\]
it follows that $f_P\leq h_P$.  By \eqref{E:hat5} and
\eqref{E:hat6}, if one considers $\phi^{-1}$, it follows that
$f_P^{-1} \leq h_P^{-1}$. Since the functions $f_P,h_P:[0,1]\to
[0,1]$ are monotone increasing we then conclude that $f_P=h_P$.
From the inequality
\begin{equation*}
\begin{gathered}
f_P( \lambda) \phi(P)=\phi(\lambda P)\leq \
\phi( \lambda (P + Q))\leq \\
\phi( \lambda P +Q) =h_P( \lambda)  \phi(P)+\phi(Q)= f_P( \lambda)
\phi(P)+\phi(Q)
\end{gathered}
\end{equation*}
we infer that
\[
0\leq \phi( \lambda(P+Q))-f_P( \lambda)  \phi(P)\leq \phi(Q).
\]
As $\phi(Q)$ is of rank one, this implies that
\begin{equation}\label{E:hat7}
\phi( \lambda (P+Q))= f_P( \lambda)  \phi(P) +k_P( \lambda)
\phi(Q)
\end{equation}
holds for some scalar $k_P( \lambda) \in [0,1]$.

With the notation $F=P+Q$ it follows from the equality \eqref{E:hat7} that
the operator $\phi(\lambda F)$ is diagonizable with respect to any
orthonormal basis in the range of the projection
$\phi(F)=\phi(P)+\phi(Q)$. We obtain that $\phi(\lambda F)$ is a
constant multiple of $\phi(F)$. This gives us that we have $f_P=k_P$
and hence
\[
\phi( \lambda (P+Q))= f_P( \lambda)  \phi(P) +f_P( \lambda)
\phi(Q).
\]
Let $R$ be a rank-one projection whose range is included in the subspace
generated by $\rng P$ and $\rng Q$. Then we have
\[
f_R(\lambda)\phi(R)=
\phi(\lambda R)\leq
\phi(\lambda F)=
f_P(\lambda)(\phi(P)+\phi(Q))=
f_P(\lambda)\phi(F).
\]
This gives us that
\begin{equation}\label{E:hat30}
f_R \leq f_P
\end{equation}
whenever $P,R$ are rank-one projections. It is then trivial
that there is in fact equality in \eqref{E:hat30} and this proves that
$f_P$ does not depend on $P$. Denote by $f$ this
common function.

Our next claim is to show that $f$ satisfies a functional equation
of the form \eqref{E:FE}. Fix mutually orthogonal rank-one projections
$P,Q$ on $H$. Pick $\mu \in ]0,1[$ and let $E=\mu P +Q$. Take any
rank-one projection $R$ on $H$ whose range is contained in the
subspace generated by the ranges of $P$ and $Q$ and which is
neither equal nor orthogonal to $P$. Using the formula
\eqref{E:hat3}, one can easily verify that
\[
\lambda (E,R)= \frac{\mu}{\mu +(1-\mu)\tr PR}.
\]
By the definition of $\lambda(E,R)$ and the order preserving
property of $\phi$ it is clear that
\begin{equation*}
\begin{gathered}
f(\lambda(E,R))= \sup\{ f( \lambda)  \, : \,  \lambda R\leq E\}=
\sup\{ f( \lambda)  \, : \, \phi( \lambda R)\leq \phi(E)\}=\\
\sup\{ f( \lambda)  \, : \, f( \lambda)  \phi(R)\leq \phi(E)\}=
\lambda(\phi(E),\phi(R)).
\end{gathered}
\end{equation*}
Since $\phi(E)=\phi(\mu P+Q)=f(\mu) \phi(P)+\phi(Q)$, we obtain
the equality
\[
f\biggl( \frac{\mu}{\mu+(1-\mu)\tr PR}\biggr)=
\frac{f(\mu)}{f(\mu) +(1-f(\mu))\tr \phi(P)\phi(R)}.
\]
As the quantities $\tr PR$ and $\tr\phi(P)\phi(R)$ do not depend
on $\mu$, it follows from this equality that $\tr \phi(P)\phi(R)$
can be uniquely expressed as a function of $\tr PR$. Denoting
$g(\tr PR)=\tr \phi(P)\phi(R)$, we get a bijective function $g
:]0,1[ \to ]0,1[$ for which
\[
f\biggl( \frac{\mu}{\mu+(1-\mu)\nu}\biggr)= \frac{f(\mu)}{f(\mu)
+(1-f(\mu))g(\nu)} \qquad (\mu,\nu \in ]0,1[).
\]
This gives us the desired functional equation for $f$ and $g$. By
Proposition~\ref{P:hat8} we obtain that there are positive real numbers
$a,b,c$ such that
\[
  f(x)=\frac{x^c}{x^c+a(1-x)^c} \qquad(x\in]0,1[)
\]
and
\[
  g(y)=b y^c\qquad(y\in]0,1[).
\]

But our function $g$ has the additional property that
$g(1-x)=1-g(x)$ $(x\in ]0,1[)$. In fact, this follows from the equality
\begin{equation*}
\begin{gathered}
g(\tr PR)+ g(1-\tr PR)=g(\tr PR)+g(\tr QR)= \\
\tr \phi(P)\phi(R)+\tr \phi(Q)\phi(R)=1.
\end{gathered}
\end{equation*}
One can easily deduce that we necessarily have $b=1, c=1$,
i.e., $g$ is the identity on $]0,1[$. This further implies that
our function $f$ is of the form
\[
f(x)=\frac{x}{x+a (1-x)}=\frac{x}{x(1-a)+a} \qquad (x\in
]0,1[).
\]
Because of continuity, the above equality holds also on the whole interval $[0,1]$.
Hence, we have that $f$ is of the form $f=f_p$ where $p=1-a$.

Since the function $g$ above is the identity, we have
\[
\tr PQ=\tr \phi(P)\phi(Q)
\]
for all rank-one projections $P,Q$ on $H$. Hence, using Wigner's
celebrated theorem on symmetry transformations (sometimes called
uni\-ta\-ry-antiunitary theorem) we obtain that
there exists an either unitary or antiunitary operator $U$ on $H$
such that
\[
\phi(P)=UPU^*
\]
holds for every rank-one projection $P$ on $H$. Consider the
transformation
\[
\psi: A\longmapsto f^{-1} (U^* \phi(A)U)
\]
on $E(H)$. It is not hard to see that this map is a
bijection of $E(H)$ which preserves the order (as well as zero
product) in both directions and it has the additional property
that it fixes the so-called weak atoms, that is, the effects of the form
$\lambda P$ where $\lambda \in [0,1]$ and $P$ is a rank-one projection.
As, according to \cite[Corollary 3]{BusGud99}, every effect is equal to the
supremum of the set of all weak atoms it majorizes, we have that $\psi$ is
the identity on $E(H)$. Transforming back, we see that
\[
\phi (A)=Uf(A) U^* \qquad (A\in E(H)).
\]
The proof is complete.
\end{proof}

\begin{proof}[Proof of Proposition~\ref{P:hat1}]
Let $\phi:E(H)\to E(H)$ be an order preserving bijection and
$\lambda ,\mu$ be a pair of nontrivial scalars for which
we have $\phi(\lambda I)=\mu I$. Keeping the notation introduced
in the first part of this section, we claim that
\begin{equation}\label{E:hat4}
f_P(\lambda)=\mu,
\end{equation}
i.e., that $\phi(\lambda P)=\mu \phi(P)$ holds for every rank-one
projection $P$. Indeed, from
\[
f_P(\lambda)\phi(P)=\phi(\lambda P)\leq \phi(\lambda I)=\mu I
\]
we deduce that $f_P(\lambda) \leq \mu$. Now, considering
$\phi^{-1}$ and $\phi(P)$ in the place of $\phi$ and $P$,
respectively, we also have $f_P^{-1}(\mu)\leq \lambda$. Since
$f_P$ is increasing, this implies $\mu\leq f_P(\lambda)$ and hence
we get \eqref{E:hat4}.

We next assert that $\phi$ preserves the orthogonality between
rank-one projections. To see this, let $P,Q$ be mutually
orthogonal rank-one projections. Denote by $P'=I-P$ the orthogonal
complement of $P$. Consider the effect $E=\lambda P+ P'$. Clearly,
we have $\lambda I\leq E\leq I$, the strength of $E$ along $P$ is
$\lambda$ and along $Q$ (which is a subprojection of $P'$) is $1$.
It follows from the order preserving property of $\phi$ and from
\eqref{E:hat4} that
\begin{itemize}
\item[--] $\mu I\leq \phi(E)\leq I$,
\item[--] the strength of $\phi(E)$ along $\phi(P)$ is $\mu$,
\item[--] the strength of $\phi(E)$ along $\phi(Q)$ is $1$.
\end{itemize}
Now, Lemma 3 in \cite{ML01d} tells us that in this case the ranges
of $\phi(P)$ and $\phi(Q)$ are subspaces of the eigenspaces of
$\phi(E)$ corresponding to the eigenvalues $\mu$ and $1$,
respectively. This yields that the ranges of $\phi(P)$ and
$\phi(Q)$ are orthogonal to each other.

As every projection is equal to the supremum of the set of
all rank-one projections it
majorizes, it follows that $\phi$ preserves the orthogonality of
projections of any rank. It is also easy to verify that $\phi$
preserves the range projections of effects.
This means that if
$R$ is the range projection of $A$ (i.e., the projection onto
${\overline{\rng A}}$), then $\phi(R)$ is the range
projection of $\phi(A)$. Indeed, this preserver property follows from the
simple fact that the range projection of the effect $A$ is equal
to the infimum of the set of all projections which are greater than or
equal to $A$. It is clear that for any $A,B \in E(H)$, we have
$AB=0$ if and only if the range projections of $A$ and $B$ are
orthogonal. Using these observations we can infer that
\[
AB=0 \Longleftrightarrow \phi(A)\phi(B)=0.
\]
This completes the proof.
\end{proof}

\begin{proof}[Proof of Corollary~\ref{C:hat5}]
By \cite[Lemma 2]{ML01d} an effect is coexistent with every other
effect if and only if it is a scalar multiple of the identity.
This implies that our transformation $\phi$ maps scalar operators
to scalar operators. Hence, by Corollary~\ref{C:hat3} we infer
that up to unitary-antiunitary equivalence, $\phi$ is of the form
\[
\phi(A)=f_p(A) \qquad (A\in E(H))
\]
for some $p<1$. We claim that $p=0$, i.e., $f_p$ is the identity.
To see this, let $P,Q$ be different rank-one projections. It was
proved in \cite[Lemma 2]{ML01d} that two rank-one effects with
different ranges are coexistent if and only if their sum is an
effect. Since $\lambda P$ and $(1-\lambda)Q$ are always coexistent
(indeed, their sum is an effect), we obtain that $\phi(\lambda P)$
and $\phi((1-\lambda) Q)$ must be also coexistent. By the just mentioned
characterization we
obtain that for any $\lambda \in ]0,1[$ we have
\[
f_p(\lambda)P+f_p(1-\lambda)Q=\phi(\lambda
P)+\phi((1-\lambda)Q)\leq I.
\]
If we let $Q$ tend to $P$, we infer from this inequality that
\begin{equation}\label{E:hat8}
f_p(\lambda)+f_p(1-\lambda)\leq 1 \qquad (\lambda \in [0,1]).
\end{equation}
Since $\phi^{-1}$ has the same properties as $\phi$, it also
follows that
\[
f_p^{-1}(f_p(\lambda))+f_p^{-1}(1-f_p(\lambda))\leq 1.
\]
This further implies
\[
f_p^{-1}(1-f_p(\lambda))\leq 1-\lambda
\]
and by the monotonicity of $f_p$ we obtain
\[
1-f_p(\lambda)\leq f_p(1-\lambda).
\]
Comparing this with \eqref{E:hat8}, we see that
\[
f_p(1-\lambda)=1-f_p(\lambda)
\]
holds for every $\lambda \in [0,1]$. It is easy to show that this
implies $p=0$ which completes the proof.
\end{proof}

\begin{proof}[Proof of Corollary~\ref{C:hat6}]
We first show that $\phi$ preserves zero product in both
directions. As $\phi$ preserves the order in both directions, we
know that $\phi$ preserves the projections in both direction.
Since $\phi$ is an ortho-order automorphism, we obtain that it
also preserves the orthogonality between projections. Now, one can use
the same argument as in the last paragraph of the proof of
Proposition~\ref{P:hat1} to verify that $\phi$ preserves zero
product in both directions. Since $\phi$ preserves the
orthocomplements, we can compute
\[
\phi\biggl(\frac{1}{2}I\biggr)=
\phi\biggl(\biggl(\frac{1}{2}I\biggr)'\biggr)=
\phi\biggl(\frac{1}{2}I\biggr)'
\]
which implies that $\phi(\frac{1}{2}I)=\frac{1}{2} I$. One can
apply Corollary~\ref{C:hat4} to complete the proof.
\end{proof}

\begin{proof}[Proof of Corollary~\ref{C:hat7}]
It was proved in \cite{GudGre02b} that the order on $E(H)$ is
completely determined by the sequential product. More precisely,
\cite[Theorem 5.1]{GudGre02b} tells us that for any $A,B\in E(H)$
we have
\[
A\leq B \Longleftrightarrow \exists\, C \in E(H) \, :\, A=B\circ
C.
\]
As $\phi$ is a sequential automorphism, using this
characterization it follows that $\phi$ preserves the order in
both directions. It is easy to see that
\[
A\circ B=(\sqrt{B}\sqrt{A})^*\sqrt{B}\sqrt{A}=0 \Longleftrightarrow AB=0.
\]
Therefore, $\phi$ also preserves zero product in both
directions. By Theorem~\ref{T:hat2} we obtain that up to
unitary-antiunitary equivalence, $\phi$ is of the form
\begin{equation}\label{E:hat31}
\phi(A)=f_p(A) \qquad (A\in E(H))
\end{equation}
for some $p<1$. We assert that $p=0$, i.e., $f_p$ is the identity.
In fact, as $\phi$ is a sequential automorphism, we obtain from
\eqref{E:hat31} that $f_p$
is a multiplicative function on the unit interval. It is then
obvious that $p=0$ and the proof is complete.
\end{proof}

In conclusion, according to our promise given in the introduction,
we make some remarks on the relation
between the class of our preservers and the collections of
ortho-order automorphisms, resp.
sequential automorphisms in the general setting of von Neumann algebras. So, let
$\mathcal A$ be a von Neumann algebra of operators acting on the
Hilbert space $H$. Denote by $E(\mathcal A)$ the set of all
effects which belong to $\mathcal A$, i.e., let $E(\mathcal
A)=E(H)\cap \mathcal A$. The definitions of the order,
orthocomplementation, and sequential product are straightforward
and so are the definitions of ortho-order automorphisms and
sequential automorphisms.

First, let $\phi$ be an ortho-order automorphism of $E(\mathcal A)$.
It is easy to see that the sharp elements in $E(\mathcal A)$
(i.e., the elements $A$ for which the infimum of $A$ and $A'=I-A$
is $0$) are exactly the projections. Hence we obtain that $\phi$
preserves the projections in $E(\mathcal A)$
as well as their orthogonality in both directions.
Since, as it is well-known, the range projection of any element
of a von Neumann algebra
also belongs to the algebra, we see that $\phi$ preserves the
range projections of the elements of $E(\mathcal A)$ in the same
sense as it was mentioned in the proof of Proposition~\ref{P:hat1}.
Then one can
argue as in that proof to verify that $\phi$ preserves zero
product in both directions. So, we obtain that every ortho-order
automorphism of $E(\mathcal A)$ belongs to our class of
preservers, that is, those automorphisms preserve the order and zero product in
both directions.

Now, let $\phi:E(\mathcal A) \to E(\mathcal A)$ be a sequential
automorphism. It is not hard to prove that the above mentioned
characterization of
the order by means of the sequential product due to Gudder and
Greechie holds true also in the setting of von Neumann algebras. This means that
for any $A,B \in E(\mathcal A)$ we have
\begin{equation}\label{E:hat33}
A\leq B \Longleftrightarrow \exists\, C \in E(\mathcal A) \, :\,
A=B\circ C.
\end{equation}
In fact, if $A,B \in E(\mathcal A)$ and $A\leq B$, then by
\cite[Theorem 5.1]{GudGre02b} there is an operator $C\in E(H)$
such that $A=B\circ C$. Following the proof in \cite{GudGre02b}, one can see that
the existence of this operator $C$ is a consequence of
a well-known result of Douglas \cite{Dou66}. In fact,
in the corresponding part of the proof of Douglas' result
this $C$ was constructed. Examining the construction, it is not
hard to verify that $C$ belongs to the von Neumann algebra $\mathcal
A$, i.e., we have $C\in E(\mathcal A)$. This gives one implication
from the asserted equivalence in \eqref{E:hat33}. The other implication is trivial.
Then, just as in the proof of Corollary~\ref{C:hat7}, one can show
that $\phi$ preserves the order and zero product in both
directions.

To sum up, it has turned out that the ortho-order
automorphisms and the sequential automorphisms all belong to our new
class of preservers even in the setting of von Neumann algebras.
This seems to be a worthwhile observation as in that generality
it is quite hard to see any connection between those two kinds of
automorphisms.

\section*{Acknowledgements}
This paper was written when the author held a Humboldt Research
Fellowship. He is very grateful to the Alexander von Humboldt
Foundation for providing ideal conditions for research, and to
his host Werner Timmermann (TU Dresden, Germany) for very warm
hospitality and, what concerns the present paper, for his kind advises
that helped to improve the presentation of the manuscript.
         The author also acknowledges support from
         the Hungarian National Foundation for Scientific Research
         (OTKA), Grant No. T030082, T031995, and from
         the Ministry of Education, Hungary, Grant
         No. FKFP 0349/2000

\newpage

\end{document}